\documentstyle[12pt]{article}

\newdimen\pIR
\title{s-Step Orthomin and GMRES implemented on parallel computers } 
\author{\large A.T.\ Chronopoulos and S.K. Kim 
\thanks{UMSI 90/43R, University of Minnesota Supercomputing Institute, Minneapolis, MN. Printed 1990. A.T. Chronopoulos,  Department of Computer Science, University of Texas, San Antonio, TX 78249, USA.  Email: anthony.chronopoulos@utsa.edu
 }\\ }
\date{}
\begin{document}
\maketitle
\begin{abstract}
The Orthomin ( {\bf Omin }) \cite{Vinsome} and the  Generalized Minimal 
Residual method ({\bf GMRES }) \cite{sa:sch} are commonly used iterative 
methods for approximating the solution of nonsymmetric linear systems. 
The s-step generalizations of these methods enhance their data locality 
parallel and properties by forming $s$ simultaneous search direction vectors. 
Good data locality is the key in achieving near peak rates on memory 
hierarchical supercomputers. 
The theoretical derivation of the s-step Arnoldi and Omin has been published 
in \cite{ch:si}, \cite{ki:ch}. Here we 
derive the s-step GMRES method. We then implement s-step 
Omin and GMRES on the Cray-2 hierarchical memory supercomputer. 

\end{abstract}

\pagebreak

\section{Introduction}
Several algorithms which improve the data locality for dense linear 
algebra problems have been suggested for shared memory systems 
(e.g.see \cite{do:so} , \cite{ja:ma}, \cite{gp:sa}). These algorithms 
are based on BLAS3 (Basic Linear algebra level 3). 
BLAS3 consist of submatrix block operations and they have 
been proven highly efficient for parallel and vector computers 
One very important advantage of these algorithms over the 
standard ones is that their very low ratio of memory references over 
floating point operations. This allows efficient use of vector registers 
and local memories. It also reduces the need for frequent synchronizations 
of the processors. Linear algebra algorithms which improve 
data locality on distributed memory systems have also been studied 
(e.g.see \cite{lb:ti}, \cite{mc:vv}). 

In the area of iterative methods for solving the linear system 
\begin{equation}
  A  x = f  
\end{equation}
with nonsingular BLAS2 modules implementations consisting of one or 
more single vector 
operations have been studied in \cite{h:vo},\cite{hv:vo}. The $s$-step 
iterative methods \cite{ch:ge},\cite{ch:gep},\cite{ch:si},\cite{ki:ch} 
use BLAS3 operations. 
Other authors have considered BLAS3 approaches for 
iterative methods \cite{j:gr}, \cite{p:sa}. The s-step methods form 
independent direction vectors using repeated matrix vector 
multiplication of the coefficient matrix with a single 
direction or residual vector. This provides coarser granularity and 
increases the parallelism by computing simultaneously the 
$2s$ inner products involved in the evaluation of parameters used in 
advancing the iterations. 

The ratio of memory references  per floating point operations  
becomes of the order $1/s$ of that of the standard 
methods in operations consisting of linear combinations and dot products. 
This also holds true for matrix times vector operations with 
narrow banded matrices. Such an important class 
of matrices is the block tridiagonal matrices obtained from discretization 
of partial differential equations on regular domains. To solve 
iteratively linear systems with these coefficient matrices on a 
hierarchical memory supercomputer 
the s-step methods would require only $1/s$ of the number of main 
memory sweeps required by the standard methods. 
Also, the  $2s$ inner products required for one $s$-step are executed 
simultaneously. This reduces the need for frequent 
global communication in a parallel system and will increase delivered 
performance. 

In this article we outline the s-step Omin and Arnoldi methods. 
We then derive the s-step GMRES method. 
We compare the s-step methods with the standard methods. The implementation 
on the Cray-2 parallel vector  processor shows that the s-step 
methods are more efficient the standard ones. 
In section 2 we review the 
Omin and GMRES methods. In section 3 we present the 
s-step Omin method. In section 4 we review the s-step 
Arnoldi method and we derive the s-step GMRES method. In sections 5 and 
6 we present numerical tests and results. In section 7 we draw conclusions. 

\section{Omin and GMRES }
In this section we describe Omin and GMRES 
\cite{ee:sch}, \cite{sa:sch}, \cite{Vinsome}. Let $ x_0$ be an initial guess to 
the solution of (1) and let $r_0 = b - A x_0$ be the initial residual. 
The Omin(k) algorithm can be summarized as follows. 
\begin{description}
\item {\bf Algorithm 2.1}   {\rm Omin(k) } 
   \item Compute $r_0 $ and set $p_0 = r_0$.
   \item {\bf For } $i = 0 , 1, \ldots $ {\bf until } convergence {\bf do }
          \begin{enumerate}
             \item[1.] $a_i = \frac {r_i^T A p_i }{ (A p_i)^T A p_i}$
             \item[2.] $x_{i+1} = x_i + a_i p_i$
             \item[3.] $r_{i+1} = r_i - a_i A p_i$
             \item[4.] $b_j^i = \frac {(A r_{i+1})^T A p_j}
                   { (A p_j)^T A p_j }$ for $j_i \leq j \leq i$
             \item[5.] $p_{i+1} = r_{i+1} - \sum_{j = j_i}^i b_j^{(i)} p_j$
             \item[6.] $Ap_{i+1} =Ar_{i+1} - \sum_{j= j_i}^i b_j^{(i)} Ap_j$
          \end{enumerate}
   \item {\bf EndFor }
\end{description}

In this algorithm $ j_i $ = $min(0, i-k+1)$ for Omin(k). Fixing 
$j_i  =  0 $  yields the Generalized Conjugate Residual algorithm ({\bf GCR})
\cite{ee:sch} which is equivalent to GMRES. 
For Computational cost we count only vector operations as inner products, 
vector updates and matrix vector products ({\bf Mv }). 
For Omin(k) each iteration (with $  k \leq i-1 $ ) needs: $k+2$ inner 
products, $2k+2$ vector updates and 1 Mv. 

The GMRES method \cite{sa:sch} is based on the Arnoldi procedure 
for computing an $l_2 -$orthonormal basis 
{$q_1 , q_2 , \ldots ,  q_j$} of the Krylov subspace ${\bf K_m}$ = 
${\rm span}  \{ q_1, Aq_1, \ldots A^{m-1} q_1  \} $
If $Q_j$ is the $n \times j$ matrix whose columns are the $l_2 -$
orthonormal basis 
$\{ q_1 , q_2 , \ldots ,  q_j \}$, then $H_j = {Q_j}^T AQ_j $,
is the upper $j \times j$ Hessenberg matrix whose entries are the scalars 
$h_{i,l}$ generated by the Arnoldi method. GMRES consists of 
an Arnoldi procedure and an error minimization step. 
We next present one cycle of the restarted GMRES(m) method. The 
norm of the residual is monitored for the convergence check.

\begin{description}
\item {\rm {\bf Algorithm 2.2} GMRES($m$)}
   \item  Compute $r_0$ and set $q_1 = \frac{ r_0}{ \| r_0 \|} $
   \item  {\bf For } $j = 1, \ldots , m-1$
          \begin{enumerate}
             \item[1.]  $h_{i,j} = q_i^T A q_j, \; 1 \leq i \leq j$
             \item[2.]  $\hat{q}_{j+1} = A q_j
                                            - \sum_{i=1}^{j} h_{i,j} q_i$
             \item[3.]  $h_{j+1,j} = \| \hat{q}_{j+1} \|_2$
             \item[4.]  $q_{j+1} = \hat{q}_{j+1} / \| \hat{q}_{j+1} \|_2$
          \end{enumerate}
   \item  {\bf EndFor }
\end{description}
Form the approximate solution: $x_m  =  Q_m y_m $,  
where $y_m $ minimizes 
\begin{equation}
 J(y)  =  \| \beta e_1 - G_m y \|  \;\;\;\; e_1 = [ 1 , \ldots , 0 ]^T  . 
\end{equation} 
The matrix $ G_m $ is 
the same as $ H_m $ except for an additional row whose only nonzero 
element is $ h_{m+1,m} $ in the $ ( m+1 , m )$ position. 
Minimizing the error functional $m$-dimensional $J(y)$ is equivalent to 
solving: 
\begin{equation}
 \min_{x \in x_0 + {\bf K}_m} \| b - Ax \|_2 
\end{equation}
where ${\bf K}_m = {\rm span} \{ r_0, Ar_0, \ldots A^{m-1} r_0 \}$ 
is the Krylov subspace of dimension $m$. The linear 
least squares problem (2) is solved by use of the $QR$ method. More details 
can be found in \cite{sa:sch}. For GMRES each iteration needs: $i+2$ inner 
products, $i+1$ vector updates and 1 Mv. 

\section{s-Step Omin}
The s-step Minimal Residual ({\bf MR}) method is a simple steepest descent 
method which computes the following sequence of solution approximations 
\[ x_{i+1}  =  x_i  +  a_i^1 r_i  +   \ldots   +  a_i^s A^{s-1} r_i , \] 
where $ a_i^j $ to minimizes $\| r_{i+1} \| $ over 
the affine Krylov subspace 
\[ \{ x_i  +  \sum_{j=0}^{s-1} a_j A^j r_
i    :  a_j \;\; {\rm scalars \;\;}  {\rm and \;\;} r_i  =  f  -  A x_i \} \]
This method is theoretically equivalent to GMRES. Unlike GMRES 
s-MR is not be stable for large $s$ because of loss of orthogonality 
of the direction vectors used. 

The s-step MR is used to obtain s-step generalizations for GCR and 
Omin(k). The details can be found in \cite{ch:si}. To achieve this 
we form the s directions $ \{ \; r_i , \ldots , A^{s-1} r_i \;\} $ and 
simultaneously $A^T A$ -orthogonalize to 
$k$ preceding blocks of direction vectors 
$ \{ \;\; [ p_j^1 , \ldots , p_j^s ] \;\; \}_{ j=j_i }^{j=i} $. 
The norm of the residual $ \| r_{i+1} \|_2 $ is minimized 
simultaneously in all $s$ new directions in order to obtain $x_{i+1}$. 
The following notation (in BLAS3) facilitates the description of the 
algorithm. 
\begin{itemize}
   \item Set $P_i = [ p_i^1 , \ldots , p_i^s ]$ 
   \item Set $R_i = [ r_i , A r_i ,..., A^{s-1} r_i ] $.
   \item Let $W_i $ be
         $W_i = [ \; (A p_i^j)^T A p_i^l \; ], \;\; 1 \leq j, \; l \leq s.$
   \item Let $\underline{a}_i, \underline{m}_i $ be the vectors 
         $\underline{a}_i = [ a_i^1 , \ldots , a_i^s ]^T$ 
         and $\underline{m}_i = [ r_i^T A p_i^1 , \ldots , r_i^T A p_i^s ]^T.$
   \item For $l = 1 , \ldots , s$ and $ j = j_i , \ldots , i$
         let $\underline{c}_j^l, {\underline{b}_j^l }$ be the  
         vectors  $\underline{b}_j^l = [ b_j^{(l,1)} , \ldots , b_j^{(l,s)}]^T$
         and
\[\underline{c}_j^l = [ A^{(l+1)T} r_{i+1} A p_j^1 ,
                   \ldots , A^{(l+1)T} r_{i+1}, A p_j^s ) ]^T .\]
\end{itemize}
Using BLAS3 operations we summarize s-step Omin(k) in the following algorithm. 
\begin{description}
\item {\bf Algorithm 3.1} {\rm $s$-step Orthomin($k$)  }
   \item Compute $R_0$ and set $P_0 = R_0 $. 
   \item {\bf For } $i = 0, 1, \ldots $ {\bf until } convergence {\bf do }
      \begin{enumerate}
         \item Compute $\underline{m}_i , W_i$.
         \item (Scalar1) Decompose $W_i$ and solve
               $W_i \underline{a}_i = \underline{b}_i$.
         \item $x_{i+1} = x_i + P_i \underline{a}_i$
         \item $r_{i+1} = r_i - A P_i \underline{a}_i$
         \item Compute $R_i$.
         \item Compute $\underline{c}_j^i$, for $j = j_i ,..., i$
         \item (Scalar2) Solve $W_j \underline{b}_j^l = - \underline{c}_j^l$,
               for $j = j_i , \ldots , i$ and $l = 1 , \ldots , s$.
         \item $P_{i+1} = R_{i+1} + \sum_{j = j_i}^i P_j 
                  \underline{b}_j^l $.
         \item $AP_{i+1} = AR_{i+1} + \sum_{j = j_i}^i AP_j 
                  \underline{b}_j^l $.
      \end{enumerate}
   \item {\bf EndFor }
\end{description}
The value of the index $ j_i $ is $ {\rm min} ( 0 , i-k+1) $. 
Fixing $j_i = 0 $ yields the s-step GCR method. 
For $s=1$ we obtain the standard Omin(k) and GCR 
methods. It is proved in \cite{ch:si} that s-step Omin(1) 
coincides with s-step GCR if $A $ is symmetric or skew-symmetric. 
$s$ must be kept small (not greater than five) for numerical 
stability reasons \cite{ch:ge}.  However if 
$R_i$ is made into an $A^T A - 
$orthogonal set (by use of modified Gramm-Schmidt) after it is computed a  larger $s$ can be chosen. 
In this case it can be easily shown that the linear systems in Scalar1, 2 
are diagonal. We do not go into details here because we have not yet 
implemented this approach. 

In \cite{ch:si} it is proved that s-step Omin(k) converges  
for nonsymmetric definite matrices and for a class of indefinite matrices. 
In fact since each iteration contains an s-MR iteration it converges  
for the same class of matrices as GMRES(s).  

We give the vector work and storage for Omin(k) and s-step Omin(k) in Table 3.1. Storage includes the matrix $A$ and the vectors: 
\[ x, \; r, \; AR, \; \{  P_j  \} _{j=0}^{j = i+1} , \; \{ AP_j \}_{j=0}^ {j = i+1} . \] 
More details can be found in \cite{ch:si}.

\begin{center}

{\bf Table 3.1: Ops/storage for j-iter of Omin(k) and s-Omin(k) 
 }
\begin{tabular}{||c c c ||} 
\hline
Vector Ops& Omin((k) &s-Omin(k) \\ [0.5ex]
\hline\hline
Dotprod &min($[(j+1)+2],[k+2]$)&min($[(j+1)s^2+s(s+1)/2],$ \\
& & $[ks^2+s(s+1)/2]$) \\
\hline
Matvec & 1  &	s \\
\hline
Vect Update&min($[2(j+1)+1],[2k+2]$)&min($[2(j+1)s^2+s], $\\
& & $[2ks^2+s(s+1)/2]$) 	\\
\hline
Storage & matr A+(2k+2)vects	&matr A+(2ks+s +1)vects\\[1ex]
\hline
\end{tabular}
\end{center}

\section{s-step GMRES }
The s-step Arnoldi method has been derived in \cite{ki:ch}. 
Firstly we give an outline the s-step Arnoldi method and then we derive 
residual error minimization step which yields the s-step GMRES method. 

Let us denote by $k$ the iteration number in the $s$-step Arnoldi method.
Given the vectors $\{ \;\;{v_k}^1 , {v_k}^2 , \ldots , {v_k}^s \;\;\} $
(each of dimension $N$) we use $\bar{V}_k$ to denote the matrix of 
[${v_k}^1 , {v_k}^2 , \ldots , {v_k}^s$]. 
Initially we start with a vector ${v_1}^1$ and compute 
${v_1}^2  =  A {v_1}^1 \;\;$, ... ,$\;\; {v_1}^s  =  A^{s-1} 
{v_1}^1$. One way to obtain an $s$-step Arnoldi algorithm is to use 
these $s$ linearly independent vectors and generate a sequence of 
block matrices $\bar{V}_1 , \ldots $. To form $\bar{V}_2$ we compute  
${v_2}^1  =  A {v_1}^s$,${v_2}^2  =  A {v_2}^1 \;\;$, ... ,$\;\; {v_2}^s  =  A^{s-1} {v_2}^1$.
Then we orthogonalize $\bar{V }_2$ against $\bar{V }_1$.
Inductively we form $\bar{V }_k$ for $k>1$. 
The subspaces $\bar{V }_1$,$\bar{V }_2$, \ldots $\bar{V }_k$ are mutually orthogonal, but the vectors 
{${v_k}^1$,$ {v_k}^2$, \ldots ,$ {v_k}^s$} are not orthogonal, that is,
${\bar{V }_k}^T \bar{V }_k$ is not a diagonal matrix. 

We next summarize the s-step Arnoldi algorithm using BLAS3 operations. 
 
\begin{description}
\item {\bf Algorithm 4.1} {\rm s-step Arnoldi }
   \item  Select $v_1^1 $
   \item  Compute $\bar{V} _1 = [ $ ${v_1}^1$,$ {v_1}^2  = A{v_1}^1$, \ldots ,$ {v_1}^s  = {A^{s-1}} {v_1}^1 ] $  
    \item  For $k = 1, \ldots, m/s$
          \begin{enumerate} 
             \item      Call Scalar1 
             \item      Compute $ v_{k+1}^1  =  A v_k^s  -  \sum_{i=1}^k \bar{V} _i [{\bf  {\bar{h }}_{ik}^1  }] $. 
             \item      Compute $ {v_{k+1}}^2  = A{v_{k+1}}^1$, \ldots ,$ {v_{k+1}}^s  = {A^{s-1}} {v_{k+1}}^1  $  
             \item      Compute $(A^i {v_k}^1 , {v_l}^j ) \;\; {\rm for}  1 \leq i, j \leq s \; {\rm and} 1 \leq l \leq k-1      $
             \item      Compute $(A^i {v_k}^1 , A^j {v_k}^1 ) \;\; {\rm for}  0 \leq i \leq s-1  \;  {\rm and} i \leq j  \leq s $
             \item      Call Scalar2 
             \item      Compute $\bar{V} _{k+1}  =  [ v_{k+1}^1 , \ldots , v_{k+1}^s ]  -  \sum_{i=1}^k \bar{V} _i [ 0 ,  {\bf t_{ik}^1} , \ldots , {\bf t_{ik}^{s-1}} ]$  
\end{enumerate}
    \item  {\bf EndFor }
\end{description}
\begin{description}
\item {\bf Scalar1: } Compute and decompose $W_i  =  \bar{V} _i^T \bar{V}  _i $. 
\item solve $W_i {\bf {\bar{h }_{ik}  }^q}  =  {\bf {b_{ik}}^q} $ for $q  =  1, \ldots ,s$, where 
\item $ {\bf {b_{ik}}^1}  =  { [ ( {v_i}^1 ,  A{v_{k+1}}^1 ) ,  \ldots  , ( {v_i}^s ,  A{v_{k+1}}^1 ) ] }^T ,  \ldots  , $ 
\item $    {\bf {b_{ik}}^{s-1}}  =   [ ( {v_i}^1 ,  {A^{s-1}}{v_{k+1}}^1 ) ,  \ldots  , ( {v_i}^s ,  {A^{s-1}}{v_{k+1}}^1 ) ] ^T $ 
\item {\bf Scalar2: } Solve $W_i  {\bf t_{ik}  ^q}  =   {\bf {b_{ik}}^q} $ for $1 \leq q \leq s-1$, where  
\item $ [  {\bf c_{ik}^1 } , \ldots , { \bf c_{ik}^s } ]   =  \bar{V} _i^T A \bar{V} _k $ 
\end{description}

It was proven in \cite{ch:si} that the inner products computed in steps 5. and 6. and 
scalar work can be used in evaluating the coefficients and right hand 
sides of the linear systems in Scalar 1, 2. The parameters 
$ {\bf {\bar{h }_{ik}  }^q}$ computed in Scalar 1 are not only 
important computing $ v_{k+1}^1 $ but also they are entries of the 
upper Hessenberg matrix $\bar{H}_k $ of the s-step Arnoldi method. 
The following matrix equality holds but it is not explicitly computed 
except for the vector $ v_{k+1}^1 $ (in step 2. of algorithm 4.1):
\begin{equation}
                 A \bar{V} _k  =   \sum_{j=1}^k  \bar{V} _j [ {\bf \bar{h}_{jk}^1},  \ldots  ,  {\bf \bar{h}_{jk}^s } ]  + v_{k+1}^1 e_{sk}^T
\end{equation} 
The upper Hessenberg matrix $\bar{H}_k $ is obtained 
from matrix  $\bar{H}_{k-1}$ by adding the block column 
$ [ {\bf \bar{h}_{jk}^1},  \ldots  ,  {\bf \bar{h}_{jk}^s } ]$ for 
$j= 1, \ldots ,k$ 
down to the diagonal plus an $s \times s$ subdiagonal block which has the only 
one nonzero $\| v_{k+1}^1 \| _2 ^2 $ at the position $(s(k-1)+1, sk)$. 

To introduce an $s$-step GMRES method we use the basis generated by 
the $s$-step Arnoldi method. 
After $k$ iterations of $s$-step Arnoldi method 
we have ${v_{k+1}}^1$ and a $(k+1)s \times ks$ matrix $\bar{G }_k$ 
generated by the method. $\bar{G}_k$ is the same as $\bar{H }_k$
except for an additional row whose only nonzero element is at the $( s(k-1)+1 ,  s k)$ position.

Let $V_k  = $ $[ \bar{V}_1 , \ldots , \bar{V}_k ] $ and 
$U_k  = $ $[ V_k , {v_{k+1}}^1] $ 
then the matrix $\bar{G }_k$ satisfies the important relation:
\begin{equation}
A V_k = \bar{G}_k U_k
\end{equation}
To derive the s-step GMRES we must solve the least squares problem:
\begin{equation}
\min_{ z \in K_j }  \| f - A [ x_0  + z] \| =  \min_{ z \in K_j } 
\|  r_0  - Az \|.
\end{equation}
If we set $z = {V_k} \bar{y }$, we can view the norm to be minimized as 
the following function of $\bar{y }$:
\begin{equation}
J(y) = \|  {v_1}^1  -  A V_k \bar{y } \|
\end{equation}
where we have let ${v_1}^1  =   r_0 $ for convenience. Using euation (5) 
\begin{equation}
J(y) = \| U_k [ e_1  -  \bar{G }_k \bar{y } ] \|.
\end{equation}
Here the vector $e_1  =  [ 1 , \ldots , 0 ]^T $.  

Let $D = W_k^T W_k$, then
$D = diag( {\bar{V }_1}^T \bar{V }_1, \ldots , \bar{{V }_k}^T \bar{V }_k, {{v_{k+1}}^1}^T {v_{k+1}}^1)$.
Using Cholesky factorization $D=L^T L$, we obtain
\begin{equation}
J(y) = \|  [ \beta e_1  -  L \bar{G}_k \bar{y } ] \|.
\end{equation}
where $\beta = \| r_0 \| = \|{v_1}^1 \| $ because $ L e_1  =  
\| u_1^1 \| e_1 $. Hence the solution of the least squares 
problem (6) is given by 
\begin{equation}
x_k  =  x_0  +  V_k \bar{y }_k
\end{equation}
where $\bar{y }_k$ minimizes the functional $J(y)$. 
 
Now we describe the restarting s-step GMRES algorithm.

\begin{description}
\item {\bf Algorithm 4.2} {\rm s-step GMRES }
          \begin{enumerate}  
   \item  Compute $r_0$ and set $v_1^1 = \frac{ r_0}{ \| r_0 \|} $ 
   \item  Compute the s-step Arnoldi vectors $V_1 \; \ldots \; V_m$ 
   \item  Form the approximate solution: $x_m  =  x_0  +  V_m \bar{y }_m$ 
   \item  where $y_m $ minimizes $J(y) = \|  [ \beta e_1  -  L \bar{G}_k \bar{y } ] \|$. 
  \item   Restart: Compute $r_m  =  f  -  A x_m$ and stop if $\| r_m \| < \epsilon $ else set $ x_0  =  x_m $ and $r_0 = r_m $ and go to 1. 
          \end{enumerate}
\end{description}

The advantages of $s$-step GMRES compared to standard GMRES on parallel 
computers come from the fact the matrix vector operations, inner products 
and linear combinations corresponding to s consecutive steps of the 
standard GMRES(m) are grouped together for simultaneous execution. 
Scalar1 and Scalar2 and the minimization of $J(y)$ are scalar computations 
of dimension $s$.  
 
{\bf  Remark 4.1:  } Assume that the degree of the minimal polynomial 
of $r_0 $ is greater than $ ms $. Let GMRES(ms) and 
s-step GMRES(m) start with the same $x_0 $. Then the iterate 
$x_{i(ms)} $ is the same for the s-step GMRES(m) and GMRES(ms). 

This result follows from the equivalence of these methods to the $ms$-step MR method. 
We compare the computational work and storage of the $s$-step GMRES method to 
the standard one. We give the vector work for the $s$-step and standard GMRES in Table 4.1. We present the storage and  vector operations of 1 cycle of the standard GMRES($sm$) compared to 1 cycle of of the s-step GMRES(m). The details of deriving the formulas are in 
\cite{sa:sch} and \cite{ki:ch1}.
\begin{center}
{\bf Table 4.1:Vector Ops of GMRES(sm) vs s-GMRES(m)} 

\begin{tabular}{||c c c ||} 
\hline
Vector Ops& GMRES(sm) &s-GMRES(m) \\ [0.5ex]
\hline\hline
Dotprod & $ms +[sm(ms+1)/2]$ & $ [m(m-1)s^2]/2+[s(s+1)/2+s]$  \\
\hline
Matvec & (ms+1)  &	s(m+1) \\
\hline
 Vect update& $([m^2s^2+ms]/2 +2ms$ & $ m(m+1)s^2$	\\
\hline
  Storage &Matr A+ $(ms+1) vect $	&$Matr A+([s(m+1)m]/2+m)vect$	\\[1ex]
\hline
\end{tabular}
\end{center}

\section{Numerical Tests}

We have discretized a boundary value problems in partial differential 
equations on a square region by the method of 
finite differences.

{\bf  Problem:  }
 
\begin{equation}
  -{(b(x,y) u_x )}_x - {(c(x,y) u_y )}_y  + ({d(x,y) u)}_x  + 
 {(e(x,y) u)_y }  + f(x,y) u  = g(x,y)  ,
\end{equation}
\begin{equation}
 \Omega  = (0,1) \times (0,1)
\end{equation}
where
$b(x,y) = e^{-xy} , c(x,y) =  e^{xy}  , d(x,y) =  \beta (x+y) $
 
$e(x,y) =  \gamma   (x+y) , f(x,y) = \frac{1}{(1 + x y)} , $
 
$u(x,y) =  x {e^{xy}} sin( \pi y) sin( \pi y) , $

This problem is a standard elliptic test problem 
which can be found in \cite{sa:sch} and the right hand side function is 
constructed so that the analytic solution is known. 
The right hand side function $g(x,y)$ is obtained by applying the 
differential operator to $u(x,y)$. Dirichlet boundary conditions are 
imposed. By controlling  $\gamma $ and $\beta$, we could change the degree 
of nonsymmetry of the discretization matrix. We set $ \gamma  = 50.0 ,  
\beta  = 1.0$. We have used the five point difference operator for the 
Laplacian, central difference for the first derivative. 
For initial value, we have chosen x(i) = 0.05*mod(i,50). 

For a nonsingular matrix $K$ then the transformed system 
\begin{equation}
 [ A K ] K^{-1} x = f  
\end{equation}
is a right preconditioned form of the 
original linear system. We use right preconditioning, since it minimizes 
the residual norm rather than minimizing the norm of $r_i  $, 
where $ r_i $ is the i-th residual vector.
We use ILU(0) preconditioning in vectorizable form \cite{h:vo}. 

\section{Results}
   
We used the Cray-2 supercomputer at the Minnesota Supercomputer Institute. 
The Cray-2 is a four-processor (MIMD) supercomputer. 
All processors have equal access 
to a central memory of 512 Megawords. Each Cray-2 processor has  
8 vector registers (each 64 words long) and has
data access through a single path between its vector registers and main memory.
Each processor has 16 Kwords of local memory with no direct
path to central memory but with a separate data path between local memory
and its vector registers, and 
the six parallel vector pipelines: common memory to vector register 
(LOAD/STORE), vector register to local memory (LOAD/STORE) , floating point 
ADD/SUBTRACT, MULTPLY/DIVIDE, Integer ADD/SUBTRACT and LOGICAL  elines.  
It is possible to design assembly language 
kernels which exhibit a performance commensurate with the 4.2 nanosecond cycle
time of the Cray-2 if the computations allow it.
This means that a rate
of 459 Megaflops is possible on one processor if all arithmetic  pipelines can be kept busy. 

The maximum performance of the Cray-2 for specific applications comes from
data movement minimization, good vectorization and division into   rocessing tasks.
Because of single paths between vector register and central or local memory on
the Cray-2 system, memory transfers constitute a severe bottleneck 
for achieving maximum performance. 
Therefore, minimization of data movement results in faster execution times.

The termination criterion used  was $ \| r_i \|^{1/2}  <  
10^{-6} $.  The number of grid point in the $x$ and $y$ directions 
taken is $nx=  64 ,  128 ,  192   , 256 $. 
 
The termination criterion used  was $ || r sub i || sup 1/2 ~<~ 
10 sup -6 $.  The number of grid point in the $x$ and $y$ directions 
taken is $nx= 32, ~ 64 ,~ 128 ,~ 192 ~ , 256 $. 
 
 The selection of $s$ and $k$ in 
s-Omink) minimizes the number of iterations for each problem. 
Since in each iteration of s-GMRES(m) there is an overhead of one matrix 
vector multiplication (see table 4.1) $s$ must be chosen  as large as 
possible. 
So we chose first $m$ in GMRES(m) and then choose $s$ and $m bar $ (in 
s-GMRES( $m bar$ ) so that $m= m bar ~ s $. 
In our tests the number of iterations in GMRES(10) equals 
(in almost all cases) the number of iterations 5-GMRES(2). 
This is expected from remark 4.1. This is not always true 
with Orthomin((k+1)s-1) and s-step Orthomin(k). The residual 
error is minimized on the same number of independent vectors. However, these 
vectors do not generate the same affine subspace. 
Actually, s-Omin(k) may converge for indefinite problems for which 
Orthomin(k) fails [3].

Tables 6.1-6.4 contain the no. of iterations and total execution times for convergence of the 
methods. The performance gain is about 1.5 for 2-Orthomin(2) and 
1.3 for 2-GMRES(5). The 2-GMRES(5) has one additional matrix 
vector operation and since the preconditioner has low Megaflop rate 
it offsets the gains made from the other types of operations. 
This accounts for the low rate of 2-GMRES(5) in the preconditioned case. 
In terms of 
programmer optimizations we have unrolled the linear combinations 
as single GAXPY operations. Unrolling loops for BLAS3 operations did not 
lead to a faster rate. 

\begin{center}
{\bf Table 6.1:s-Omin(k) vs s-GMRES(m): Iter's ILU(0) precond. }
 \begin{tabular}{||c c c c c||} 
 \hline
 dimension & s=1,k=4 & s=2,k=2 & s=1,m=10 & s=2,m=5  \\ [0.5ex] 
 \hline\hline
 64 & 193& 98& 30& 30  \\
 \hline
 128 &335&167&60&59  \\
 \hline
 192 &509&252&93&95  \\
 \hline
 256 &660&340&146&146   \\ [1ex] 
 \hline
\end{tabular}
\end{center}

\begin{center}
{\bf Table 6.2:Omin(4) vs 2-Omin(2) times(sec) Cray-2 p-PEs  }
 \begin{tabular}{||c c c c c||} 
 \hline
 dimension & p=1& p=4 & p=1 & p=4  \\ [0.5ex] 
  \hline\hline
 64 & 0.289&0.602& 0.229& 0.503  \\
 \hline
 128 &1.947&1.503&1.495&1.138  \\
 \hline
 192 &7.169&3.113&5.175&2.168  \\
 \hline
 256 &17.523&6.873&12.374&4.827   \\ [1ex] 
 \hline
\end{tabular}
\end{center}

\begin{center}
{\bf Table 6.3:GMRES(10) vs 2-GMRES(5) times(sec) Cray-2 p-PEs  }
 \begin{tabular}{||c c c c c||} 
 \hline
 dimension & p=1& p=4 & p=1 & p=4  \\ [0.5ex] 
  \hline\hline
 64 & 0.534&1.069& 0.381& 0.560  \\
 \hline
 128 &3.495&2.392&2.782&1.471  \\
 \hline
 192 &12.253&5.831&9.183&3.855 \\
 \hline
 256 &33.741&12.378&25.978&9.051   \\ [1ex] 
 \hline
\end{tabular}
\end{center}

\begin{center}
{\bf Table 6.4:s-Omin(k) vs s-GMRES(m):time(sec), p=1, Cray-2 }
 \begin{tabular}{||c c c c c||} 
 \hline
dimension & Omin(4) & 2-Omin(2) & GMRES(10) & 2-GMRES(5)  \\ [0.5ex] 
  \hline\hline
 64 & 0.137&0.079& 0.223& 0.205  \\
 \hline
 128 &0.899&0.613&1.487&1.326 \\
 \hline
 192 &3.216&1.827&4.863&3.893 \\
 \hline
 256 &7.372&4.519&10.917&9.074   \\ [1ex] 
 \hline
\end{tabular}
\end{center}

\section{Conclusions}

We reviewed the standard and s-step Omin and Arnoldi methods 
(which are derived in \cite{ch:si}, \cite{ki:ch} ). We have derived the 
s-step GMRES method. 
We then used these methods to solve nonsymmetric systems 
arising from the finite difference discretization of partial differential 
equations. The s-step methods showed similar convergence properties as the 
standard methods with gains in execution time. The gains would have 
been much higher with use of BLAS3 modules to utilize the local 
memory of the Cray-2 processor. If the matrix vector multiplication is 
more costly then the gains from the s-step methods are reduced. 
However one could program computers with local memory 
to make matrix vector multiplication with block tridiagonal matrices 
very fast. On the Cray-2 
this has to be done in Cray Assembly language (CAL). We are currently 
implementing (in CAL) a matrix vector multiplication module which 
utilizes efficiently the local memory of the Cray-2. This is expected 
to increase the speed of the matrix vector operations for this type 
of block tridiagonal matrices.

\end{document}